\RequirePackage{plautopatch}
\documentclass[11pt,reqno]{amsart}

\usepackage[a4paper,margin=28mm]{geometry}
\usepackage[T1]{fontenc}
\usepackage{lmodern}
\usepackage{amsmath,amssymb,amsthm,mathtools}
\usepackage{array,booktabs}
\usepackage{microtype,graphicx,scalerel,xcolor}
\usepackage[pdfusetitle]{hyperref}
\usepackage{url}

\hypersetup{
    colorlinks=true,
    pdftitle={L-systems on the large Witt design},
    pdfauthor={Toshihiro Shimizu and Norihide Tokushige},
}

\usepackage{tikz}
\usetikzlibrary{positioning}
\usetikzlibrary{svg.path}
\definecolor{orcid_color}{HTML}{A6CE39}

\hypersetup{pdfborder={0 0 0}} 
\DeclareRobustCommand{\orcidicon}{%
	\raisebox{.2mm}{\scalerel*{%
	\begin{tikzpicture}[xscale=1,yscale=-1,transform shape]
	\filldraw[color=orcid_color] svg {M256,128c0,70.7-57.3,128-128,128C57.3,256,0,198.7,0,128C0,57.3,57.3,0,128,0C198.7,0,256,57.3,256,128z};
	\filldraw[color=white] svg {M86.3,186.2H70.9V79.1h15.4v48.4V186.2z} svg {M108.9,79.1h41.6c39.6,0,57,28.3,57,53.6c0,27.5-21.5,53.6-56.8,53.6h-41.8V79.1z M124.3,172.4h24.5
		c34.9,0,42.9-26.5,42.9-39.7c0-21.5-13.7-39.7-43.7-39.7h-23.7V172.4z} svg {M88.7,56.8c0,5.5-4.5,10.1-10.1,10.1c-5.6,0-10.1-4.6-10.1-10.1c0-5.6,4.5-10.1,10.1-10.1
		C84.2,46.7,88.7,51.3,88.7,56.8z};
	\end{tikzpicture}}{|}}%
}
\newcommand{\orcid}[1]{\href{https://orcid.org/#1}{\orcidicon}}

\def\F{\mathbb F}
\def\K{\mathbb K}
\def\R{\mathbb R}

\def\P{\mathbb P}
\def\Eq{\mathrm{Eq}}
\DeclareMathOperator{\supp}{supp}
\DeclareMathOperator{\wt}{wt}
\DeclareMathOperator{\rk}{rk}
\DeclareMathOperator{\rank}{rank}
\DeclareMathOperator{\cl}{cl}
\DeclareMathOperator{\Ing}{Ing}
\DeclareMathOperator{\IS}{IS}

\newcommand{\BB}{\mathcal B}
\newcommand{\CC}{\mathcal C}
\newcommand{\FF}{\mathcal F}
\newcommand{\GG}{\mathcal G}
\newcommand{\II}{\mathcal I}
\newcommand{\JJ}{\mathcal J}
\newcommand{\MM}{\mathcal M}
\newcommand{\PP}{\mathcal P}

\newtheorem{theorem}{Theorem}
\newtheorem{lemma}{Lemma}
\newtheorem{proposition}[lemma]{Proposition}
\newtheorem{corollary}[theorem]{Corollary}

\title{$L$-systems on the large Witt design} 
\author[Toshihiro Shimizu]{Toshihiro Shimizu\,\orcid{0009-0009-1379-0148}}
\address{Fujitsu Research, Fujitsu Limited, Kawasaki 211-8588, Japan}
\email{t.shimizu\_2@fujitsu.com}
\author[Norihide Tokushige]{Norihide Tokushige\,\orcid{0000-0002-9487-7545}}
\address{College of Education, University of the Ryukyus, Nishihara  903-0213, Japan}
\email{hide@cs.u-ryukyu.ac.jp}
\urladdr{https://n-tokushige.github.io/x/}
\subjclass{05D05, 05B05, 05B35}
\keywords{Witt design, Steiner system, information inequality, matroid}

\begin{document}

\begin{abstract}
Let $L=\{0,1,2,3,4,8\}$, associated with the large Witt design $S(5,8,24)$, and let
$m(n,24,L)$ be the maximum size of a family of $24$-element subsets of an $n$-element set
whose distinct members have intersection sizes in $L$. We prove that there are positive
constants $c,C,n_0$ such that, for every $n\geq n_0$,
\[
cn^{316/63}\leq m(n,24,L)\leq Cn^{169/29}.
\]
We also obtain similar results for the systems associated with $S(4,7,23)$, $S(3,6,22)$,
and $S(3,4,14)$. The upper bounds follow from a non-Shannon information inequality and
disprove conjectures of Kupavskii and Terekhov. The lower bounds follow from an algebraic
construction for arbitrary Steiner systems.
\end{abstract}
\maketitle

\section{Introduction}
We study the maximum size of $L$-systems related to Steiner systems including the large
Witt design. For upper bounds, we use a non-Shannon information inequality introduced by
Zhang and Yeung in \cite{ZY}, and further developed by Dougherty, Freiling, and Zeger in
\cite{DFZ} and by Csirmaz in \cite{Cs}. For lower bounds, we give an algebraic
construction with an explicit exponent, using the dimension formula for generic
intersections due to Brakensiek, Gopi, and Makam \cite[Theorem~1.16]{BGM2}.

For a finite set $X$, let $2^X$ denote the power set of $X$, and $\binom Xi$ denote the
set of $i$-element subsets of $X$, where $i$ is an integer with $0\leq i\leq |X|$. For a
positive integer $n$, write $[n]:=\{1,2,\ldots,n\}$. Let $k$ be a positive integer, and
let $L\subset\{0,1,\ldots,k-1\}$. We say that $\FF\subset\binom{[n]}k$ is an $L$-system,
or a $(k,L)$-system if $|F\cap F'|\in L$ for all distinct $F,F'\in\FF$. We write
$m(n,k,L)$ for its maximum possible size. Deza, Erd\H{o}s, and Frankl \cite{DEF}
introduced this problem and proved that, for fixed $k,L$ and sufficiently large $n$,
\begin{align*}
m(n,k,L)\leq\prod_{l\in L}\frac{n-l}{k-l}.  
\end{align*}
Some basic results concerning $m(n,k,L)$ are summarized in \cite{FOT,T} by Frankl, Ota,
and Tokushige.

In this paper, we consider the situation $n\to\infty$ for fixed $k,L$. If there is a
constant $\alpha$ such that $m(n,k,L)=\Theta(n^\alpha)$, then we call $\alpha=\alpha(k,L)$
the exponent of the $(k,L)$-system. Such an exponent exists and is rational for all
$k\le12$ and $L\subset\{0,1,\ldots,k-1\}$ \cite{FOT,T}, but its existence for arbitrary
$k,L$ remains open. Recently, Kupavskii and Terekhov \cite[Corollary~1.3]{KT} showed that
\[
 \gamma(k,L):=\lim_{n\to\infty}\frac{\log m(n,k,L)}{\log n}
\]
exists for all $(k,L)$. Thus $m(n,k,L)=n^{\gamma(k,L)+o(1)}$, and
$\alpha(k,L)=\gamma(k,L)$ whenever $\alpha(k,L)$ exists. The quantity denoted by $\alpha$
in \cite{KT} is our $\gamma$. One of our main results is the following: if
$L=\{0,1,2,3,4,8\}$, then
\[
5+\frac1{63}\leq\gamma(24,L)\leq 6-\frac5{29}. 
\]
This $L$-system is related to the large Witt design, as we will explain shortly.

Let $t,b,k$ be integers with $k>b>t\geq 2$. We say that a family $\BB\subset\binom{[k]}b$
is a Steiner system $S(t,b,k)$ if every $t$-subset of $[k]$ lies in exactly one member of
$\BB$. Members of $\BB$ are called blocks, and $B(T)$ denotes the unique block containing
$T\in\binom{[k]}t$. We also call $\BB$ a design (or a block design) of $S(t,b,k)$. Two
designs with the same parameters need not be isomorphic. The Steiner systems $S(5,6,12)$
and $S(5,8,24)$ are each unique up to isomorphism. They are called the small and large
Witt designs, respectively; see, e.g., \cite{MW-S}.

For a Steiner system $S=S(t,b,k)$, set
\[
 L(S):=\{0,1,\ldots,t-1,b\}.
\]
It is known from \cite{FOT,T} that if $k\leq 12$, then $\alpha(k,L(S))=|L(S)|=t+1$. In
particular, it is shown in \cite{T} that, for the small Witt design $S=S(5,6,12)$,
\[
 \alpha(12,\{0,1,2,3,4,6\})=6.
\]
The large Witt system behaves differently. Our first theorem shows that
$\gamma(k,L(S))<t+1$ if
\begin{align}\label{eq:1}
S\in\{S(5,8,24),\,S(4,7,23),\,S(3,6,22),\,S(3,4,14)\}.
\end{align}

\begin{theorem}\label{thm1}
For each row of the following table, there is a constant $C=C(k,L)>0$ such that
$m(n,k,L)\leq Cn^\beta$ for every $n\geq k$.
\[
\begin{array}{ccc}
\toprule
k&L&\beta\\
\midrule
24&\{0,1,2,3,4,8\}&169/29=6-5/29\\
23&\{0,1,2,3,7\}&140/29=5-5/29\\
22&\{0,1,2,6\}&111/29=4-5/29\\
14&\{0,1,2,4\}&111/29=4-5/29\\
\bottomrule
\end{array}
\]
In particular, $\gamma(k,L)\leq\beta$ in all four cases.
\end{theorem}

For a given $(k,L(S))$-system $\FF$, we can extract a subfamily $\FF^*$ with
$|\FF^*|\geq c|\FF|$ and a closed intersection structure by F\"uredi's
theorem~\cite{Furedi}. In essential cases, the intersection structure contains a Steiner
design $\BB$ of $S(t,b,k)$. Then, $\FF^*$ can be viewed as a $k$-partite family (which we
call a $\BB$-agreement code) with some intersection property. The design defines an
associated matroid rank function. The rank and the code entropies are then related, and a
non-Shannon inequality concerning entropies gives an upper bound for $|\FF^*|$.

Theorem~\ref{thm1} refutes Conjectures 1.5 and 11.3 of \cite{KT}, and the associated
matroids refute Conjecture 11.5. We explain these consequences in subsection~
\ref{sec:KT conjectures}.

\medskip
Our second theorem gives an explicit lower bound for $m(n,k,L(S))$. In $S(t,b,k)$, the
number of blocks containing a fixed $u$-subset, where $0\le u\le t$, is given by
\begin{align}\label{eq:2}
\lambda_u:=\binom{k-u}{t-u}\bigg/\binom{b-u}{t-u}. 
\end{align}

\begin{theorem}\label{thm2}
Suppose that $k>b>t\geq2$ and that $S=S(t,b,k)$ exists. Then there are constants $c>0$ and
$n_0$ such that, for all $n\geq n_0$,
\[
m(n,k,L(S))\geq c n^{\,t+(t-1)/(\lambda_1-1)}.
\]
Consequently,
\[
\gamma(k,L(S))\geq t+\frac{t-1}{\lambda_1-1}.
\]
\end{theorem}

\begin{corollary}
The lower bounds in the following table hold.
\[
\begin{array}{ccc}
\toprule
k&L(S)&\text{lower bound for }\gamma(k,L(S))\\
\midrule
24&\{0,1,2,3,4,8\}&316/63=5+1/63\\
23&\{0,1,2,3,7\}&307/76=4+3/76\\
22&\{0,1,2,6\}&31/10=3+1/10\\
14&\{0,1,2,4\}&77/25=3+2/25\\
\bottomrule
\end{array}
\]
\end{corollary}

R\"odl and Tengan \cite{RT} have already established a lower bound
$\gamma(k,L(S))\geq t+\epsilon$ for some $\epsilon=\epsilon(t,b,k)>0$. They first obtain a
seed family using a near-perfect matching theorem with a probabilistic proof, and then
enlarge it by a product construction. Our construction is algebraic and gives the
displayed value of $\epsilon$. It is related in spirit to the construction using
generating matrices in \cite{T}. Our seed family is obtained by a finite-field realization
from a generic matrix based on a result by Brakensiek, Gopi, and Makam \cite{BGM2}. Then
we enlarge the seed family by a finite-field extension.

\section{Upper bounds}
In this section, all random variables are discrete and finite.

\subsection{Non-Shannon information inequalities}
We recall the notation for entropy and mutual information; see
\cite[Chapter~2]{CT}. All logarithms are the natural logarithms.
For a finite random variable $X$ with support $\mathcal X:=\{x:\P[X=x]>0\}$, define the
entropy of $X$ by
\[
H(X):=-\sum_{x\in\mathcal X}\P[X=x]\log \P[X=x].
\]
Joint entropy is defined in the same way for a tuple of random variables. We use the
convention $0\log0=0$. If $Y$ has support $\mathcal Y$, define the conditional entropy of
$X$ given $Y$ by
\[
H(X\mid Y):=\sum_{y\in\mathcal Y}\P[Y=y]H(X\mid Y=y).
\]
The identity $\P[X=x,Y=y]=\P[Y=y]\P[X=x\mid Y=y]$ gives the chain rule
\[
H(X,Y)=H(Y)+H(X\mid Y).
\]

The following elementary bounds are proved in Appendix~\ref{appendix:A}.
\begin{lemma}\label{lem:1}
If $X$ takes at most $q$ values, then $0\leq H(X\mid Y)\leq H(X)\leq\log q$.
\end{lemma}

Mutual information and conditional mutual information are defined by
\begin{align}
I(X;Y)&:=H(X)+H(Y)-H(X,Y),\nonumber\\
I(X;Y\mid Z)&:=H(X,Z)+H(Y,Z)-H(Z)-H(X,Y,Z).\label{eq:3}
\end{align}
By expanding the definitions and applying Jensen's inequality to the concavity of the
logarithm, we obtain the basic information inequality
\begin{align}\label{eq:4}
I(X;Y\mid Z)\geq0. 
\end{align}
Entropy, conditional entropy, mutual information, and conditional mutual information are
nonnegative. Such inequalities asserting these nonnegativities are called the basic
inequalities. Informally, information inequalities implied by the basic inequalities are
called Shannon-type inequalities; see \cite[Chapters 13--14]{Y} for a formal definition.
There are, however, valid information inequalities that are not Shannon-type. The first
such non-Shannon-type inequality was found by Zhang and Yeung \cite{ZY}:
\[
 -I(A;B)+2I(A;B\mid C)+I(A;C\mid B)+I(B;C\mid A)+I(A;B\mid D)+I(C;D)\geq 0.
\]
We use the above inequality and its relatives to prove Theorem~\ref{thm1}. To this end,
define the Ingleton expression by
\[
\Ing(A,B,C,D):=-I(A;B)+I(A;B\mid C)+I(A;B\mid D)+I(C;D),
\]
and, for $i_2\geq i_1>0$, $i_3,i_4\geq0$, let
\begin{align}
\MoveEqLeft{f(i_1,i_2,i_3,i_4;A,B,C,D)}\nonumber\\
&:=  i_1\Ing(A,B,C,D)+(i_2-i_1)I(A;B\mid C)+i_3I(A;C\mid B)+i_4I(B;C\mid A)\label{eq:5}\\
&= -i_1I(A;B)+i_2I(A;B\mid C)+i_3I(A;C\mid B)+i_4I(B;C\mid A)\nonumber\\
&\qquad +i_1I(A;B\mid D)+i_1I(C;D).\nonumber
\end{align}
Fix $A,B,C,D$ in this order, and write $f(i_1,i_2,i_3,i_4)$ for
$f(i_1,i_2,i_3,i_4;A,B,C,D)$.

\begin{lemma}[\cite{ZY,DFZ,Cs}]\label{lem:2}
We have $f(i_1,i_2,i_3,i_4)\geq0$ if
\[
(i_1,i_2,i_3,i_4)\in\{(1,2,1,1),\,(8,17,9,5),\,(5,9,5,5)\}. 
\]
\end{lemma}

The first quadruple gives the Zhang--Yeung inequality \cite{ZY}, in the equivalent form
given in \cite[Theorem~3]{DFZ}. The second is entry 82) in the table in Section~VIII,
page~19, of \cite{DFZ} due to Dougherty--Freiling--Zeger. The third is the final
unnumbered inequality in \cite[Section~6.5]{Cs}, due to Csirmaz. Here we only include a
quick proof of the case $(1,2,1,1)$ for readers' convenience. (The other two cases can be
shown in a similar manner, but the procedure is more involved.) For the proof we need the
following version of the Copy Lemma \cite[Lemma~2]{DFZ}, see also \cite[Section~6]{Cs}.
Juxtaposition below denotes a tuple, e.g., $AB$ means $(A,B)$.

\begin{lemma}[Dougherty--Freiling--Zeger]\label{lem:3}
Let $A,B,C,D$ be jointly distributed random variables. Then there is another random
variable $R$, jointly distributed with $A,B,C,D$ with the following properties.
\begin{itemize}
 \item The marginal distributions of $(A,B,C)$ and $(A,B,R)$ are the same with $R$
replacing $C$. In particular,
\begin{align}\label{eq:6}
 H(XR)=H(XC)\text{ where } X\in\{\emptyset,A,B,AB\}.
\end{align}
\item The variables $(C,D)$ and $R$ are independent over $(A,B)$, that is,
\begin{align}\label{eq:7}
 I(CD;R\mid AB) = 0.
\end{align}
\end{itemize}
\end{lemma}

\begin{proof}[Proof of the Zhang--Yeung inequality]
We show that $f(1,2,1,1)\geq 0$, where
\[
f(1,2,1,1)=-I(A;B)+2I(A;B\mid C)+I(A;C\mid B)+I(B;C\mid A)+I(A;B\mid D)+I(C;D).
\]
Expansion of $f(1,2,1,1)$ by the chain rule reduces to
\begin{align}\label{eq:8}
\begin{split}
&-2(H(A)+H(B))-H(C)+3(H(AB)+H(BC)+H(AC))\\
&\qquad +H(AD)+H(BD)-H(CD)-4H(ABC)-H(ABD). 
\end{split}
\end{align}
Let $R$ be the random variable in Lemma~\ref{lem:3}. Then the chain rule gives
\begin{align}\label{eq:9}
I(D;R\mid AB)=0, \quad I(C;R\mid AB)=0. 
\end{align}
Indeed, by \eqref{eq:7}, $0=I(CD;R\mid AB)=I(D;R\mid AB)+I(C;R\mid ABD)$. Since both terms
on the right are nonnegative, $I(D;R\mid AB)=0$. Similarly, $I(C;R\mid AB)=0$. Now put
\begin{align*}
\Delta&:=
I(C;R\mid A)+I(C;R\mid B)+I(A;B\mid CR)+I(D;R\mid A)\\
&\qquad +I(D;R\mid B)+I(A;B\mid DR)+I(C;D\mid R)+I(AB;R\mid CD).
\end{align*}
Then $\Delta\geq 0$ because all terms in the RHS are nonnegative. Expansion of $\Delta$ by
the chain rule with
\eqref{eq:6}, \eqref{eq:7}, and \eqref{eq:9} reduces to 
\eqref{eq:8}. Thus, $f(1,2,1,1)\geq 0$ follows.
\end{proof}

\subsection{Rank function and entropies}
Let $\BB$ be a design of a Steiner system $S(t,b,k)$. Let $\Omega=[k]$ and let $q\geq2$ be
an integer. We call $[q]$ an alphabet, $\CC\subset[q]^\Omega$ a code, and
$x=(x_i)_{i\in\Omega}$ a word. For $x,y\in\CC$, let $\Eq(x,y):=\{i\in\Omega:x_i=y_i\}$ be
the set of agreement coordinates. Let
\begin{align}\label{eq:10}
 \II_\BB:=\bigcup_{i=0}^{t-1}\binom\Omega{i}\cup\BB. 
\end{align}
We say that $\CC\subset[q]^\Omega$ is a $\BB$-agreement code if $\Eq(x,y)\in\II_\BB$ for
all distinct $x,y\in\CC$.

Suppose that $\CC$ is a $\BB$-agreement code, where $\BB$ is a design of $S(t,b,k)$.
Define a function $r_\BB:2^\Omega\to\{0,1,\ldots,t+1\}$ by
\begin{align*}
 r_\BB(U):=
\begin{cases}
 |U| & \text{ if }|U|<t,\\
 t & \text{ if $|U|\geq t$ and there exists $B\in\BB$ such that $U\subset B$},\\
 t+1 & \text{ otherwise.}
\end{cases}
\end{align*}
This is the rank function of the matroid whose independent sets are
\[
\JJ=\bigcup_{i=0}^t\binom\Omega{i}\ \cup\ 
\left\{U\in\binom\Omega{t+1}:U\not\subset B
\text{ for every }B\in\BB\right\}.
\]
We do not need this observation to prove Theorem~\ref{thm1}, but we only mention it in for
Subsection~\ref{sec:KT conjectures}.

If $r_\BB(U)=t$, then we write $B(U)$ for the unique $B\in\BB$ with $U\subset B$. Indeed,
for every $T\in\binom Ut$, there is exactly one $B\in\BB$ such that $T\subset B$, because
$\BB$ is a design of $S(t,b,k)$. In particular, if $|U|=t$, then $r_\BB(U)=t$. For
$U\subset\Omega$, let $x_U:=(x_i)_{i\in U}$. Define an equivalence relation $x\sim_U y$ on
$\CC$ by $x_U=y_U$. If $r_\BB(U)=t$, then $U$ and $B(U)$ induce the same equivalence
classes. Define a projection $\pi_U:\CC\to[q]^U$ by $\pi_U(x)=x_U$.

Choose one word $X=(X_i)_{i\in\Omega}$ uniformly from $\CC$. Throughout the upper-bound
argument, every $X_U$ is the coordinate projection of this same random word. In
particular, its distribution is the induced marginal distribution and need not be uniform
on $\pi_U(\CC)$.

Define a normalized entropy
\begin{align*}
 h(U):=\frac{H(X_U)}{\log q}, 
\end{align*} 
and let
\begin{align*}
 R:=\frac{H(X)}{\log q}=\log_q|\CC|.
\end{align*}

\begin{lemma}\label{lem:4}
Let $\CC$ be a $\BB$-agreement code, and let $U,V\subset\Omega$.
\begin{itemize}
 \item[(i)] If $r_\BB(U)=t$, then $h(U)=h(B(U))$, and $h(U)\leq t$.
 \item[(ii)] If $r_\BB(U)=t+1$, then $h(U)=R$.
 \item[(iii)] If $r_\BB(V)=r_\BB(U)+1=t$ and $U\subset V$, then $0\leq h(V)-h(U)\leq 1$.
\end{itemize}
\end{lemma}

\begin{proof}
(i) Since $U$ and $B(U)$ induce the same equivalence classes, they give the same entropy.
That is, $H(X_U)=H(X_{B(U)})$. Choose $T\in\binom{B(U)}t$. Then $H(X_T)=H(X_{B(U)})$.
Thus, $h(U)=h(B(U))=h(T)$. Since $X_T$ takes at most $q^t$ values, Lemma~\ref{lem:1} gives
us that $H(X_T)\leq \log q^t$ and $h(T)\leq t$. Hence, $h(U)=h(T)\leq t$.

(ii) In this case, $\pi_U$ is injective, and so $H(X_U)=H(X)=R\log q$.

(iii) Let $B:=B(V)$ and choose $x\in B\setminus U$. By $r_\BB(U)=t-1$ and the definition
of $r_\BB$, we have $|U|=t-1$. Then $U$ and $T:=U\cup\{x\}$ determine the same $B$. By
(i), $h(T)=h(B)=h(V)$. The chain rule gives
$H(X_x\mid X_U)=H(X_x,X_U)-H(X_U)=H(X_T)-H(X_U)$. Since $X_x$ takes at most $q$ values, it
follows from Lemma~\ref{lem:1} that $0\leq H(X_x\mid X_U)\leq \log q$ and
$0\leq h(T)-h(U)\leq 1$. This together with $h(T)=h(V)$ gives the required bound.
\end{proof}

For brevity, we write $AB$ for $A\cup B$, $ABC$ for $A\cup B\cup C$, etc.

\begin{proposition}\label{prop:5}
Let $\Omega=[k]$ and let $\BB$ be a design of $S(t,b,k)$. Suppose that
$A,B,C,D\in\binom\Omega{t-1}$ satisfy
\begin{align*}
& r_\BB(AB)= r_\BB(AC)= r_\BB(AD)= r_\BB(BC)= r_\BB(BD)=t,\\
& r_\BB(CD)= r_\BB(ABC)= r_\BB(ABD)=t+1.
\end{align*}
Assume also that the information inequality $f(i_1,i_2,i_3,i_4;Y_1,Y_2,Y_3,Y_4)\geq0$
holds for every four jointly distributed finite random variables $Y_1,Y_2,Y_3,Y_4$. Then,
every $\BB$-agreement code $\CC\subset[q]^\Omega$ satisfies
$|\CC|\leq q^{(t+1)-\epsilon}$, where
\begin{align}\label{eq:11}
 \epsilon=\frac{i_1}{2i_1+i_2+i_3+i_4}. 
\end{align}
\end{proposition}

\begin{proof}
Let $\delta:=t+1-R$. The tuple $(X_A,X_B)$ and $X_{A\cup B}$ determine one another, even
when $A$ and $B$ overlap. Consequently, $H(X_A,X_B)=H(X_{A\cup B})$. With abuse of
notation, write $A,B,C,D$ for $X_A,X_B,X_C,X_D$. Also, put $I^*(\cdot):=I(\cdot)/\log q$
and $\Ing^*(\cdot):=\Ing(\cdot)/\log q$.

We claim that $0\le I^*(A;B\mid C)\leq\delta$. The lower bound follows from \eqref{eq:4}.
For the upper bound, we first expand $I^*(A;B\mid C)$ by \eqref{eq:3}, and apply
Lemma~\ref{lem:4} with the rank assumptions. Then,
\[
 I^*(A;B\mid C)=\left(h(AC)-h(C)\right)+h(BC)-h(ABC)
 \leq 1+t-R=\delta.
\]
In this way, we have
\begin{equation}\label{eq:12}
 0\le I^*(A;B\mid C),\ I^*(A;C\mid B),\ I^*(B;C\mid A)\le\delta.
\end{equation}
Similarly,
\begin{align}
 \Ing^*(A,B,C,D)
 &=\left(h(AB)+h(AD)+h(BC)\right)-\left(h(ABC)+h(ABD)+h(CD)\right)\notag\\
 &\qquad+\bigl(h(AC)-h(A)\bigr)+\bigl(h(BD)-h(B)\bigr)\notag\\
 &\leq 3t-3R+1+1=3\delta-1.\label{eq:13}
\end{align}
Let $f^*(i_1,i_2,i_3,i_4):=f(i_1,i_2,i_3,i_4)/\log q$. Substituting \eqref{eq:12} and
\eqref{eq:13} into \eqref{eq:5} yields
\[
 0\leq f^*(i_1,i_2,i_3,i_4)\leq i_1(3\delta-1)+\left((i_2-i_1)+i_3+i_4\right)\delta
 =(2i_1+i_2+i_3+i_4)\delta-i_1.
\]
Thus $\delta\geq i_1/(2i_1+i_2+i_3+i_4)=\epsilon$, or equivalently, $R\leq (t+1)-\epsilon$, 
as required.
\end{proof}

We record the value of $\epsilon=\epsilon(i_1,i_2,i_3,i_4)$ corresponding to Lemma~\ref{lem:2}:
\[
\begin{array}{cc}
\toprule
(i_1,i_2,i_3,i_4)&\epsilon\\
\midrule
(1,2,1,1)&1/6\approx0.166\\
(8,17,9,5)&8/47\approx0.170\\
(5,9,5,5)&5/29\approx0.172\\
\bottomrule
\end{array}
\]
In the next three subsections, we verify that, for $S=S(t,b,k)$ listed in \eqref{eq:1}, we can 
choose four subsets $A,B,C,D$ satisfying the assumption of Proposition~\ref{prop:5}.

\subsection{The large Witt design}\label{sec:witt}
There is, up to isomorphism, a unique design $\BB$ of the Steiner system $S(5,8,24)$,
called the large Witt design; see \cite[Chapter~20, \S5]{MW-S}. The eight-element blocks
of the large Witt design are called octads.

Let $\Omega:=[24]$. The concrete representation of $\BB\subset\binom{\Omega}8$ is given in
Appendix~\ref{appendix:B}. Let $\CC\subset[q]^\Omega$ be a $\BB$-agreement code. Let $X$
be uniform on $\CC$, and retain the common projections $X_U$, the normalized entropy $h$,
and $R$ defined above. Let
\begin{align*}
A=\{5,8,11,14\},\,
B=\{7,8,14,16\},\,
C=\{8,12,13,14\},\,
D=\{8,14,17,19\}.
\end{align*}

\begin{lemma}\label{lem:6}
We have
\begin{align*}
& r_{\BB}(AB)=r_{\BB}(AC)=r_{\BB}(AD) =r_{\BB}(BC)=r_{\BB}(BD)=5,\\
& r_{\BB}(CD)=r_{\BB}(ABC)=r_{\BB}(ABD)=6.
\end{align*}
\end{lemma}

\begin{proof}
The following six octads belong to $\BB$, as certified in Appendix~\ref{appendix:B}:
\begin{align*}
 O_{AB}&=\{1,2,5,7,8,11,14,16\},& O_{AC}&=\{3,5,8,11,12,13,14,15\},\\
 O_{AD}&=\{4,5,8,10,11,14,17,19\},& O_{BC}&=\{7,8,12,13,14,16,18,24\},\\
 O_{BD}&=\{7,8,14,16,17,19,21,23\},& O_{T}&=\{8,9,13,14,15,17,19,24\}.
\end{align*}
If $U\in\{AB,AC,AD,BC,BD\}$, then $U\subset O_U$ and $|U|=6$. Thus, $r_\BB(U)=5$.

Let $T:=\{8,13,14,17,19\}$. Then $T\subset O_T$, and so $O_T$ is the unique octad
containing $T$. Since $T\subset CD$ but $12\in CD\setminus O_T$, no octad contains $CD$.
Also $|CD|=6$. Thus, $r_\BB(CD)=6$.

Let $T':=\{5,7,8,11,14\}$. Then $T'\subset O_{AB}$, and so $O_{AB}$ is the unique octad
containing $T'$. Since $T'\subset AB$, but $12\in C\setminus O_{AB}$ and
$17\in D\setminus O_{AB}$, no octad contains either $ABC$ or $ABD$. Also $|ABC|=|ABD|=8$.
Thus $r_\BB(ABC)=r_\BB(ABD)=6$.
\end{proof}

\subsection{Derived designs}
Let $\BB$ be the large Witt design. For $E\subset\Omega:=[24]$, let
\[
 \BB_E:=\{F\setminus E:E\subset F\in\BB\}.
\]
We will use the classical uniqueness of the designs $S(4,7,23)$ and $S(3,6,22)$ up to
isomorphism; see, e.g., \cite{B}. Thus the derived designs constructed below represent
every design with their respective parameters. We will consider $S(4,7,23)$, $S(3,6,22)$,
corresponding to $E=\{8\},\{8,14\}$, respectively. Let
\[
 s=|E|,\quad t=5-s, \quad b=8-s, \quad k=24-s.
\]

For $U\in\{A,B,C,D\}$ and $E\in\{\{8\},\{8,14\}\}$, let $U_E:=U\setminus E$, e.g.,
$A_{\{8\}}=\{5,11,14\}$.

\begin{lemma}\label{lem:7}
We have
\begin{align*}
& r_{\BB_E}(A_EB_E)=r_{\BB_E}(A_EC_E)=r_{\BB_E}(A_ED_E)
 =r_{\BB_E}(B_EC_E)=r_{\BB_E}(B_ED_E)=t,\\
& r_{\BB_E}(C_ED_E)=r_{\BB_E}(A_EB_EC_E)=r_{\BB_E}(A_EB_ED_E)=t+1.
\end{align*}
\end{lemma}

\begin{proof}
The set $E$ is contained in each of $A,B,C,D$. Let $S$ be any of the eight unions
$AB,AC,AD,BC,BD,CD,ABC,ABD$, and put $S_E=S\setminus E$. For any octad $O$ we have
$S_E\subset O\setminus E\in\BB_E$ if and only if $E\subset O\in\BB$ and $S\subset O$.
Indeed, if $S\subset O$, then $E\subset S\subset O$, and the converse follows from
$S=E\cup S_E$. Since $|S_E|\geq 6-s=t+1$, the rank of $S_E$ in the derived design is $t$
precisely when $S$ is contained in an octad, and is $t+1$ otherwise. The rank assertions
now follow simultaneously from Lemma~\ref{lem:6}.
\end{proof}

\subsection{Another Steiner system: \texorpdfstring{$S(3,4,14)$}{S(3,4,14)}}
This is the case $t=3$, $b=4$, and $k=14$. Mendelsohn and Hung classified $S(3,4,14)$ in
\cite{MH}, and they showed that there are precisely four (non-isomorphic) designs
$\BB_1,\BB_2,\BB_3$, and $\BB_4$ of $S(3,4,14)$. Appendix~\ref{app:4designs} lists four
representatives and gives finite checks of their design properties and pairwise
nonisomorphism. In this description, we define $A,B,C,D$ as follows:
\begin{align}\label{eq:14}
 (A,B,C,D)=
\begin{cases}
(\{1,7\},\{2,4\},\{5,12\},\{10,11\}) &\text{for }\BB_1,\\
(\{1,3\},\{4,5\},\{8,10\},\{9,12\}) &\text{for }\BB_2,\\
(\{1,6\},\{2,3\},\{7,9\},\{11,12\}) &\text{for }\BB_3,\\
(\{1,2\},\{4,7\},\{5,14\},\{8,11\}) &\text{for }\BB_4.
\end{cases}
\end{align}
\begin{lemma}\label{lem:8}
Let $\BB\in\{\BB_1,\BB_2,\BB_3,\BB_4\}$. Then, we have
\begin{align*}
& r_{\BB}(AB)=r_{\BB}(AC)=r_{\BB}(AD) =r_{\BB}(BC)=r_{\BB}(BD)=3,\\
& r_{\BB}(CD)=r_{\BB}(ABC)=r_{\BB}(ABD)=4.
\end{align*}
\end{lemma}

\begin{proof}
In each row of \eqref{eq:14}, the four pairs are disjoint. Each of $AB,AC,AD,BC,BD$ is a
listed block and therefore has rank $t=3$, whereas $CD$ is not a listed block and has rank
$t+1=4$, see Appendix~\ref{app:4designs}. The unions $ABC$ and $ABD$ have size $6>b$, so
their ranks are also $4$.
\end{proof}

\subsection{Proof of Theorem~\ref{thm1}}
We use the following form of F\"uredi's structure theorem~\cite[Theorem $1'$ and Remark
1]{Furedi}; see also~\cite[Theorem 16.2]{FT}. A family $\FF\subset 2^{[n]}$ is $k$-partite
if its ground set $[n]$ has a partition $V_1\sqcup\cdots\sqcup V_k$ such that each member
meets every part in one point. For such a family let $\pi:[n]\to[k]$ be a projection such
that $\pi(v)=i$ when $v\in V_i$. The family is called homogeneous if
\[
 \{\{\pi(v):v\in F\cap F'\}:F'\in\FF\setminus\{F\}\}
\]
is independent of $F\in\FF$. The common family is its intersection structure, denoted by
$\IS(\FF)$. Let $L\subset\{0,1,\ldots,k-1\}$. We say that $\II\subset 2^{[k]}$ is a closed
$L$-system if $I\cap I'\in \II$ and $|I\cap I'|\in L$ for all (not necessarily distinct)
$I,I'\in\II$.

\begin{theorem}[F\"uredi]\label{thm:Furedi}
For fixed $k\geq2$ and $L\subset\{0,\ldots,k-1\}$, there is $c=c(k,L)>0$ such that every
nonempty $(k,L)$-system $\FF$ contains a homogeneous $k$-partite subfamily $\FF^*$ with
$|\FF^*|\geq c|\FF|$ and with $\IS(\FF^*)$ a closed $L$-system.
\end{theorem}

Using this structure theorem, we can obtain an upper bound for $m(n,k,L)$ from an upper
bound for the corresponding extracted codes as follows.

\begin{lemma}\label{lem:9}
Let $L=\{0,1,\ldots,t-1,b\}$ and $\beta\geq t$. Suppose that for every design $\BB$ of
$S(t,b,k)$, every integer $q\geq 2$, and every $\BB$-agreement code $\CC\subset[q]^k$, we
have $|\CC|\leq q^\beta$. Then, for every $n\geq k$, we have $m(n,k,L)\leq Cn^\beta$,
where $C$ is a constant depending on $k,L$ only.
\end{lemma}

\begin{proof}
Let $\FF\subset\binom{[n]}k$ be a given $L$-system. By Theorem~\ref{thm:Furedi}, choose a
$k$-partite subfamily $\FF^*\subset\FF$ with $|\FF^*|\geq c|\FF|$ and
$[n]=V_1\cup\cdots\cup V_k$. Let $\II:=\mathrm{IS}(\FF^*)$ be the intersection structure.
It suffices to show that $|\FF^*|\leq n^\beta$ so that we can take $C=1/c$.

First suppose that there is $T\in\binom{[k]}t$ such that $T\not\subset I$ for all
$I\in\II$. Then the map $F\mapsto F\cap\left(\bigcup_{i\in T}V_i\right)$ is injective.
Indeed, two distinct members with the same image would have an intersection pattern in
$\II$ containing $T$. Thus, $|\FF^*|\leq\prod_{i\in T}|V_i|\leq n^t\leq n^\beta$.

Next suppose that every $T\in\binom{[k]}t$ is contained in some $I\in\II$. Since $\II$ is
a closed $L$-system, if $T\subset I\in\II$, then $|I|=b$. Moreover, if there are distinct
$I_1,I_2\in\II$ with $T\subset I_i$ ($i=1,2$), then $T\subset I_1\cap I_2$. This yields
$t\leq |I_1\cap I_2|<b$ and $|I_1\cap I_2|\not\in L$, a contradiction. Hence, for every
$T\in\binom{[k]}t$ there is exactly one $I\in\II\cap\binom{[k]}b$ such that $T\subset I$.
This means that $\BB:=\II\cap\binom{[k]}b$ is a design of $S(t,b,k)$. We can view $\FF^*$
as a $\BB$-agreement code (by renaming vertices). Applying the hypothesis with $q=n$ gives
$|\FF^*|\leq n^\beta$.
\end{proof}

\begin{proof}[Completion of the proof of Theorem~\ref{thm1}]
The designs with parameters $S(5,8,24)$, $S(4,7,23)$, and $S(3,6,22)$ are unique up to
isomorphism \cite{MW-S, B}, and the four designs in Appendix~\ref{app:4designs} exhaust
$S(3,4,14)$ by \cite{MH}. Let $\BB$ be a design of these Steiner systems $S(t,b,k)$.
Relabeling coordinates preserves the size of the $\BB$-agreement code $\CC\subset[q]^k$.
Consequently, Proposition~\ref{prop:5}, together with Lemmas~\ref{lem:6},
\ref{lem:7}, and
\ref{lem:8}, gives the required bound 
$|\CC|\leq q^\beta$, where $\beta=(t+1)-\epsilon\geq t$. Here we use
$\epsilon=\epsilon(5,9,5,5)=5/29$ given by \eqref{eq:11} with Lemma~\ref{lem:2}. Finally,
the stated upper bounds for $m(n,k,L(S))$ follow from Lemma~\ref{lem:9}.
\end{proof}

\subsection{Conjectures by Kupavskii and Terekhov}\label{sec:KT conjectures}
A function $f:2^\Omega\to\R$ is a \emph{polymatroid function} if it is normalized,
monotone, and submodular, that is,
\[
f(\emptyset)=0,\quad f(A)\leq f(B)\ (A\subset B),\quad
f(A)+f(B)\geq f(A\cup B)+f(A\cap B).
\]
Both the entropy function $h$ and the matroid rank function $r_\BB$ are polymatroid
functions. But the polymatroid axioms alone do not imply the non-Shannon inequalities used
above.

Put $\Omega:=[k]$. For a closed $L$-system $\II\subset2^{\Omega}$, let
$\MM=\II\cup\{\Omega\}$ and define
\[
\cl_\MM(A):=\bigcap\{M\in\MM:A\subset M\},\quad
\rk(\II):=\min\{|J|:J\subset\Omega,\ J\not\subset I\text{ for every }I\in\II\}.
\]
Thus $\rk(\II)$ is the minimum size of a set $J$ with $\cl_\MM(J)=\Omega$. A polymatroid
function $f$ on $\Omega$ is called compatible with $\MM$ if $f(A)=f(\cl_\MM(A))$ for all
$A\subset\Omega$. Following \cite{KT}, define $\rho(\MM)$ to be the maximum $f(\Omega)$
where $f$ runs over all polymatroid functions compatible with $\MM$ and $f(A)\leq|A|$ for
all $A\subset\Omega$. Finally, define
\begin{align*}
\rho(k,L)&:=\max_{\II}\rho(\II\cup\{\Omega\}),\qquad
\rk(k,L):=\max_{\II}\rk(\II),
\end{align*}
where both outer maxima range over all closed $L$-systems on $\Omega$. If $J$ is of
minimum size with $\cl_\MM(J)=\Omega$, and $f$ is a polymatroid function compatible with
$\MM$, then $f(\Omega)=f(J)\leq|J|$. Hence $\rho(k,L)\leq\rk(k,L)$. Also $\rk(k,L)\le|L|$.
To see this, suppose that $r>0$ and $J=\{j_1,\ldots,j_r\}$ has minimum size among sets
contained in no member of $\II$. Choose $I_s\in\II$ containing $J\setminus\{j_s\}$. Then
$j_s\not\in I_s$, and $I_1,I_1\cap I_2,\ldots,I_1\cap\cdots\cap I_r$ have strictly
decreasing sizes in $L$, proving $\rk(\II)=r\leq|L|$. Taking the maximum over $\II$ gives
$\rk(k,L)\leq|L|$.

Now suppose $S:=S(t,b,k)$ exists and let $\BB$ be a design of $S$ and $L:=L(S)$. The
family $\II_\BB$ defined by \eqref{eq:10} is a closed $L$-system. The flats (rank maximal
subsets) of $r_\BB$ are precisely $\MM:=\II_\BB\cup\{\Omega\}$. Thus $r_\BB$ is compatible
with $\MM$ and satisfies $r_\BB(A)\leq|A|$. Moreover $r_\BB(\Omega)=t+1$. Consequently,
\[
t+1\leq\rho(k,L)\leq\rk(k,L)\leq|L|=t+1.
\]
In our notation, Conjecture~1.5 in \cite{KT} asserts $\gamma(k,L)=\rho(k,L)$, and
Conjecture~11.3 in \cite{KT} asserts $\gamma(k,L)=|L|$ whenever $\rk(k,L)=|L|$.
Theorem~\ref{thm1} refutes both conjectures.

There is also a direct contradiction to Conjecture~11.5 in \cite{KT}, which asserts that
every perfect matroid design is almost entropic (see \cite{KT} for definitions). The set
of flats $\MM$ is a perfect matroid design because if $M,M'\in\MM$ with $|M|=|M'|$, then
$r_\BB(M)=r_\BB(M')$. An almost entropic function is a limit of entropy vectors, and it
must satisfy the Zhang--Yeung inequality, in particular, \eqref{eq:8} is nonnegative.
However, for the four subsets in Lemmas~\ref{lem:6}, evaluating \eqref{eq:8} by
replacing every entropy $H(X_U)$ with $r_\BB(U)$ gives $11t-5(t-1)-6(t+1)=-1<0$ ($t=5$). 
This contradicts almost entropicity, and hence this perfect matroid design is not almost 
entropic. Similarly, the perfect matroid designs corresponding to Lemmas~\ref{lem:7} and 
\ref{lem:8} are not almost
entropic.

\section{Lower bounds}
Fix a design $\BB$ of $S(t,b,k)$, where $k>b>t\geq2$. Our construction assigns a subspace
$U_i$ to each point, and a subspace $U_B$ to each block. These subspaces have nice
conditions (Lemma~\ref{lem:10}), and cosets of the $U_i$ satisfy the desired
intersections. We first construct these subspaces generically, then specialize to a finite
field. Finally we extend that field to obtain the bound for every sufficiently large $n$.

Recall the definition of $\lambda_u$ from \eqref{eq:2}. Then, for any fixed
$U\in\binom{[k]}u$, $\lambda_u=|\{B\in\BB:U\subset B\}|$. In particular, $|\BB|=\lambda_0$
and $\lambda_t=1$. Let
\[
a:=t-1, \quad d:=\lambda_1-1, \quad D:=t\lambda_1-1. 
\]

\begin{lemma}\label{lem:10}
We have
\begin{itemize}
\item[(i)] $\lambda_1\geq\binom b{t-1}\geq t$, 
\item[(ii)] $a\lambda_u\leq D-ud$ for $1\leq u\leq t$,
\item[(iii)] $a+a\lambda_1\leq D$.
\end{itemize}
\end{lemma}

\begin{proof}
(i) Fix $B\in\BB$ and $j\not\in B$. For $T\in\binom B{t-1}$, let $C_T:=B(T\cup\{j\})$. If
$T,T'\in\binom{B}{t-1}$ with $T\neq T'$ satisfy $C_T=C_{T'}$, then $|B\cap C_T|\geq t$, a
contradiction. Thus the $\binom b{t-1}$ blocks $C_T$ are distinct. All these blocks
contain $j$, so $\lambda_1\geq\binom b{t-1}$. Since $b>t$, $\binom b{t-1}\geq t$.

(ii) For $1\leq u\leq t-1$, we have $\frac{\lambda_u}{\lambda_{u+1}}=\frac{k-u}{b-u}$, and
$\delta_u:=\lambda_u-\lambda_{u+1}=\frac{k-b}{b-u}\lambda_{u+1}>0$. Then, for
$1\leq u\leq t-2$, it follows that $\frac{\delta_u}{\delta_{u+1}}=\frac{k-u-1}{b-u}\geq 1$
and $\delta_1\geq\cdots\geq \delta_{t-1}$. By comparing the average, we have
\[
 \frac{\lambda_1-\lambda_u}{u-1}=\frac1{u-1}\sum_{i=1}^{u-1}\delta_i
\geq\frac1{t-1}\sum_{i=1}^{t-1}\delta_i=\frac{\lambda_1-1}{t-1}
\]
for $2\leq u\leq t$. Including the case $u=1$, it follows that
$\lambda_u\leq\frac{t-u}{t-1}\lambda_1+\frac{u-1}{t-1}$, and
\[
 a\lambda_u\leq(t-u)\lambda_1+u-1=(t\lambda_1-1)-u(\lambda_1-1)=D-ud.
\]

(iii) By (i), we have $D-(a+a\lambda_1)=(t-a)\lambda_1-1-a=\lambda_1-t\geq 0$.
\end{proof}

Let $\F$ be a field, and let $D,M$ be positive integers. For $DM$ algebraically
independent indeterminates $X_{r,e}$, let $\K=\F(X_{r,e}:r\in[D],\,e\in[M])$ be the
rational function field. We call $W=(X_{r,e})\in\K^{D\times M}$ a generic matrix. For
$A\subset[M]$, let $W_A\leq\K^D$ denote the span of columns indexed by $A$. If
$|A|\leq D$, then $\dim W_A=|A|$.

The following result concerning dimension of generic intersection
(\cite[Theorem~1.16]{BGM2}, the full version of \cite{BGM}) is a key for constructing
large $(k,L)$-systems related to Steiner systems.

\begin{theorem}[Brakensiek--Gopi--Makam]\label{thm:BGM} 
Given $A_1, \ldots, A_l \subset [M]$ of size at most $D$, for a $D \times M$ generic
matrix $W$, we have that
\[
\dim (W_{A_1} \cap \cdots \cap W_{A_l}) = 
\max_\PP\left\{\sum_{J\in\PP}\left|\bigcap_{i\in J}A_i\right|-\left(|\PP|-1\right)D\right\},
\]
where the maximum is over all partitions $\PP$ of $[l]$ into nonempty subsets. Note that
the result is independent of the characteristic of the underlying field.
\end{theorem}

Take the column-index set $E=\BB\times[a]$, so $M=a|\BB|$, and define
\[
 I_B=\{B\}\times[a],\quad A_i=\{(B,u)\in E:i\in B\},\quad
 U_B=W_{I_B},\quad U_i=W_{A_i}.
\]
In words, $U_B$ is the subspace spanned by $a$ columns corresponding to a block $B$, and
$U_i$ is the subspace spanned by columns corresponding to all blocks $B$ containing $i$.
Then, using $D\geq a\lambda_1\geq a$ by Lemma~\ref{lem:10}, we have $\dim U_B=a$,
$\dim U_i=a\lambda_1=D-d$, and $U_B\le U_i$ whenever $i\in B$. Moreover, we shall verify
the following two conditions.
\begin{align}
&\bigcap_{i\in T}U_i=U_{B(T)} \text{ for } T\in\tbinom{[k]}t, \label{eq:15}\\
&U_B\cap U_j=\{0\} \text{ for } B\in\BB, \, j\not\in B.\label{eq:16}
\end{align}

\begin{lemma}\label{lem:11}
The generic matrix satisfies \eqref{eq:15} and \eqref{eq:16}.
\end{lemma}

\begin{proof}
Fix $T\in\binom{[k]}t$. Then, $T\subset B(T)$ and $U_{B(T)}\leq\bigcap_{i\in T}U_i$. Thus,
for \eqref{eq:15}, it suffices to show that
$\dim(\bigcap_{i\in T}U_i)\leq\dim(U_{B(T)})=a$. Let $\emptyset\neq J\subset T$. Then,
$\bigcap_{i\in J}A_i$ is the set of columns indexed by $(B,u)\in E$ with $J\subset B$ and
$u\in[a]$. Thus,
\[
\left|\bigcap_{i\in J}A_i\right|=a\lambda_{|J|}\leq D-|J|d,
\]
by Lemma~\ref{lem:10} (ii). For every partition $\PP$ of $T$, we therefore have
\[
\sum_{J\in\PP}\left|\bigcap_{i\in J}A_i\right|-(|\PP|-1)D
\leq\sum_{J\in\PP}(D-|J|d)-(|\PP|-1)D=D-td=a.
\]
Since $|A_i|=a\lambda_1=D-d<D$, Theorem~\ref{thm:BGM} applies. Taking the maximum
over all partitions gives $\dim(\bigcap_{i\in T}U_i)\leq a$, as required.

Next, to show \eqref{eq:16}, let $B\in\BB$, $j\not\in B$. Then, $I_B\cap A_j=\emptyset$ and
$|I_B\cup A_j|=a+a\lambda_1\leq D$ by Lemma~\ref{lem:10} (iii). Since any set of at most
$D$ columns of $W$ are linearly independent,
$\dim(U_B+U_j)=a+a\lambda_1=\dim U_B+\dim U_j$, yielding \eqref{eq:16}.
\end{proof}

\begin{lemma}\label{lem:12}
Let $\F_p[X]=\F_p[X_1,\ldots,X_N]$ be a polynomial ring over a prime field $\F_p$. For a
finite number of matrices $M_1(X),\ldots, M_h(X)$ over this ring, there exists a finite
extension $\F_q/\F_p$ and a substitution $X=\xi\in\F_q^N$ such that
\[
 \rank_{\F_q}M_u(\xi)=\rank_{\F_p(X)}M_u(X)
\]
for all $1\leq u\leq h$.
\end{lemma}

\begin{proof}
Let $r_u$ be the rank of $M_u(X)$ over the rational function field, and choose a nonzero
$r_u\times r_u$ minor determinant $f_u(X)$. (If $r_u=0$, set $f_u=1$.) Then,
$f(X):=\prod_{u=1}^hf_u(X)$ is a nonzero polynomial. Choose $q=p^e$ so that
$q>\deg_{X_i}f$ for all $1\leq i\leq N$. In this case, there exists $\xi\in\F_q^N$ such
that $f(\xi)\neq 0$. We show this by induction. A nonzero univariate polynomial of degree
less than $q$ cannot vanish at every element of $\F_q$. For the induction step, regard $f$
as a polynomial in $X_N$. By induction, one of its nonzero coefficient polynomials can be
made nonzero by choosing the first $N-1$ variables. The resulting nonzero univariate
polynomial has degree less than $q$, so some choice of $X_N$ makes it nonzero.

Since all $f_u(\xi)$ are nonzero, we have $\rank_{\F_q}M_u(\xi)\geq r_u$. On the other
hand, every $(r_u+1)\times(r_u+1)$ minor determinant of $M_u(X)$ is identically zero,
which remains zero after the substitution. This shows that $\rank_{\F_q}M_u(\xi)\leq r_u$.
\end{proof}

Write $G_i=W|_{A_i}$ and $G_B=W|_{I_B}$, where the vertical bar denotes restriction to
columns. For $T=\{i_1,\ldots,i_t\}$ in a fixed order, define a $(t-1)D\times ta\lambda_1$
matrix
\[
M_T:=
 \begin{bmatrix}
 G_{i_1}&-G_{i_2}&0&\cdots&0\\
 G_{i_1}&0&-G_{i_3}&\cdots&0\\
 \vdots&\vdots&\vdots&\ddots&\vdots\\
 G_{i_1}&0&0&\cdots&-G_{i_t}
 \end{bmatrix}.
\]
All $G_i$ have full column rank, because they are $D\times a\lambda_1$ matrices and
$D>a\lambda_1$ by Lemma~\ref{lem:10}.

\begin{lemma}\label{lem:13}
The generic matrix $M_T$ satisfies $\dim\ker M_T=\dim\bigcap_{i\in T}U_i$.
\end{lemma}

\begin{proof}
For $x=(x_1,\ldots,x_t)\in\K^{a\lambda_1t}$, we have
\[
M_T\,x=(G_{i_1}x_1-G_{i_2}x_2,G_{i_1}x_1-G_{i_3}x_3,\ldots,G_{i_1}x_1-G_{i_t}x_t). 
\]
Then $M_T\,x=0$ if and only if
$G_{i_1}x_1=G_{i_2}x_2=\cdots=G_{i_t}x_t\in\bigcap_{i\in T}U_i$. Thus, the map
\begin{align}\label{eq:17}
\ker M_T\ni x\mapsto G_{i_1}x_1\in \bigcap_{i\in T}U_i  
\end{align}
is an isomorphism. Indeed, injectivity follows from the full column ranks, and
surjectivity follows by expressing any vector in the intersection in each of the column
bases.
\end{proof}

This argument applies over any field whenever the matrices $G_{i_j}$ have full column
rank, and hence also applies after the specializations below. By Lemma~\ref{lem:13} with
\eqref{eq:15},  we have
\begin{align}\label{eq:18}
\rank M_T= ta\lambda_1-\dim\ker M_T=ta\lambda_1-a=(t-1)D. 
\end{align}

Fix a base field $\F=\F_p$, and a design $\BB$ of $S(t,b,k)$. Apply Lemma~\ref{lem:12}
simultaneously to the following matrices:
\begin{align*}
G_i\,\,(i\in[k]),\quad 
G_B\,\,(B\in\BB),\quad  
[G_j\,\, G_B]\,\,(B\in\BB,j\not\in B),\quad  
M_T\,\,(T\in\tbinom{[k]}t).
\end{align*}
Their generic ranks (over $\K=\F_p(X_{r,e})$) are, respectively,
\begin{align*}
 a\lambda_1,\quad
 a,\quad
 a\lambda_1+a,\quad
 (t-1)D,
\end{align*}
where we used Lemma~\ref{lem:11} (ii) for the third, and \eqref{eq:18} for the last one.
After a suitable specialization all these ranks hold over a fixed finite field $\F_q$,
where $q=p^e$. Use the same notation for the specialized matrices and subspaces. The
isomorphism \eqref{eq:17} remains valid because each $G_i$ has full column rank. Thus it
follows that
\[
 \dim\bigcap_{i\in T}U_i=ta\lambda_1-\rank M_T=ta\lambda_1-(t-1)D=a(D+1)-aD=a.
\]
Since $U_{B(T)}\leq\bigcap_{i\in T}U_i$ and $\dim U_{B(T)}=a$, we get \eqref{eq:15}.
Finally, since $[G_j\,\,G_B]$ has full column rank, \eqref{eq:16} follows.

\medskip

Let $X:=\F_q^D$. For $i\in[k]$, let $V_i:=\{(i,x+U_i):x\in X\}$ and
$V:=\bigsqcup_{i=1}^k V_i$. Then, $|V_i|=q^d$ and $|V|=kq^d$. For $x\in X$, let
$F_x:=\{(i,x+U_i):i\in[k]\}$, and define $\FF:=\{F_x:x\in X\}$.

\begin{proposition}\label{prop:14}
Let $L=\{0,1,\ldots,t-1,b\}$. The family $\FF\subset\binom Vk$ is an $L$-system and
$|\FF|=q^D$.
\end{proposition}

\begin{proof}
Let $x,y\in X$ with $x\neq y$, and $z=x-y\neq 0$. Then $x+U_i=y+U_i$ if and only if
$z\in U_i$. Thus, letting $J(z):=\{i\in[k]:z\in U_i\}$, we have $|F_x\cap F_y|=|J(z)|$. We
claim that $|J(z)|\in L$ for all $z\neq 0$. If $|J(z)|\geq t$, choose $T\subset J(z)$ with
$|T|=t$. Then, by \eqref{eq:15}, $z\in\bigcap_{i\in T}U_i=U_{B(T)}$. To show
$B(T)\subset J(z)$, suppose that $i\in B(T)$. Since $z\in U_{B(T)}\subset U_i$, we have
$i\in J(z)$. On the other hand, if $j\not\in B(T)$, then $U_{B(T)}\cap U_j=\{0\}$ by
\eqref{eq:16}. This together with $z\neq 0$ gives us $z\not\in U_j$ and $j\not\in J(z)$.
Consequently, we have $B(T)=J(z)$ and $|J(z)|=b$. Thus, $|F_x\cap F_y|\in L$ and $\FF$ is
an $L$-system. Since $k\not\in L$, $F_x=F_y$ if and only if $x=y$. In particular,
$x\mapsto F_x$ is injective, and $|\FF|=|X|=q^D$.
\end{proof}

\begin{proof}[Proof of Theorem~\ref{thm2}]
Extend the fixed field $\F_q$ to $\F_{q^l}$, $l\geq1$, without changing the matrices.
Matrix ranks are unchanged under field extension, and the same construction with 
Proposition~\ref{prop:14} gives new families $\FF^{(l)}$ on the vertex set $V^{(l)}$, where
\[
 N_l:=|V^{(l)}|=kq^{ld}, \quad |\FF^{(l)}|=q^{lD}=(N_l/k)^{D/d}.
\]
If $N_l\leq n<N_{l+1}=q^dN_l$, add isolated vertices to obtain
\[
 m(n,k,L)\geq|\FF^{(l)}|=(N_l/k)^{D/d}
>\left(\frac{n}{kq^d}\right)^{D/d}.
\]
Taking $n_0=N_1$ and $c=(kq^d)^{-D/d}$ yields the required lower bound. Indeed,
\[
 \gamma(k,L)\geq\frac Dd=\frac{t\lambda_1-1}{\lambda_1-1}=t+\frac{t-1}{\lambda_1-1}.
\]
This completes the proof of Theorem~\ref{thm2}.
\end{proof}

\subsection*{Use of AI}
The research was done under the authors' direction using GPT-5.5, GPT-5.6, and GPT-6. The
authors have independently verified and refined all mathematical arguments, and take full
responsibility for this paper.

\subsection*{Acknowledgments}
NT was supported by JSPS KAKENHI Grant Number JP23K03201.

\appendix

\section{Inequality concerning the conditional entropy}\label{appendix:A}
Recall Jensen's inequality: If $a_i>0$, $p_i\geq0$, and $\sum_i p_i=1$, then
\[
\sum_i p_i\log a_i\leq\log\left(\sum_i p_i a_i\right). 
\]

\begin{proof}[Proof of Lemma~\ref{lem:1}]
Write $p_x:=\P[X=x]$. Suppose that $|\mathcal X|\leq q$, where $\mathcal X$ is the support
of $X$. Then,
\[
H(X)=\sum_{x\in\mathcal X}p_x\log\frac1{p_x}
\leq\log\left(\sum_{x\in\mathcal X}p_x\frac1{p_x}\right)
=\log|\mathcal X|\leq\log q.
\]
Next put $p_y=\P[Y=y]$, $p_{xy}=\P[X=x,Y=y]$, and $\mathcal S=\{(x,y):p_{xy}>0\}$. On
$\mathcal S$, both $p_x$ and $p_y$ are positive, and
\[
H(X\mid Y)=\sum_{(x,y)\in \mathcal S}p_{xy}\log\frac{p_y}{p_{xy}}.
\]
Since $p_{xy}\leq p_y$, every summand is nonnegative, and $H(X\mid Y)\geq 0$. Also,
$H(X)=-\sum_{(x,y)\in \mathcal S}p_{xy}\log p_x$. Therefore
\begin{align*}
H(X\mid Y)-H(X)
&=\sum_{(x,y)\in \mathcal S}p_{xy}\log\frac{p_xp_y}{p_{xy}}
\leq\log\left(\sum_{(x,y)\in\mathcal S}p_xp_y\right)\\
&\leq\log\left(\sum_{x\in\mathcal X}p_x\sum_{y\in\mathcal Y}p_y\right)=0.
\qedhere
\end{align*}
\end{proof}

\section{The construction of the large Witt design}\label{appendix:B}
We shall construct the large Witt design as the collection of supports of the weight-eight
words of the extended binary Golay code, see e.g., \cite[Chapter~2, \S6]{MW-S}.

Define a polynomial $p(z)\in\F_2[z]$ by
\[
 p(z):=1+z+z^5+z^6+z^7+z^9+z^{11}.
\]
For a polynomial $f(z)\in\F_2[z]$ of degree at most $11$, define $c_j\in\F_2$
($0\leq j\leq 22$) by
\[
f(z)\,p(z)=\sum_{j=0}^{22}c_jz^j, 
\]
and let $v(f):=(c_0,\ldots,c_{23})\in\F_2^{24}$, where position $i\in[24]$ contains
$c_{i-1}$ and $c_{23}:=\sum_{j=0}^{22}c_j$. Then, the extended binary Golay code is
defined by
\[
\GG:=\{v(f):f=0\text{ or }\deg f\leq11\}.
\]
For $v=(v_1,\ldots,v_{24})\in\F_2^{24}$, define its support $\supp(v):=\{i:v_i=1\}$ and
the Hamming weight $\wt(v):=|\supp(v)|$. Then,
\[
\BB:=\{\supp(v):v\in\GG,\ \wt(v)=8\}
\]
is the concrete representation of the large Witt design used in this paper.

In Section~\ref{sec:witt} we defined $O_U\in\BB$ for $U\in\{AB,AC,AD,BC,BD,T\}$. The
corresponding polynomials $f_U$ are as follows:
\begin{align*}
f_{AB}(z) &:=1+z^4,&
f_{AC}(z) &:=z^2+z^3,&
f_{AD}(z) &:=z^3+z^7,\\
f_{BC}(z) &:=z^6,&
f_{BD}(z) &:=z^6+z^{11},&
f_T(z)  &:=z^7.
\end{align*}
For example,
\begin{align*}
f_{BC}(z)\,p(z)&=z^6+z^7+z^{11}+z^{12}+z^{13}+z^{15}+z^{17},\\
v(f_{BC})&=(0,0,0,0,0,0,1,1,0,0,0,1,1,1,0,1,0,1,0,0,0,0,0,1),\\ 
\supp(v(f_{BC}))&=\{7,8,12,13,14,16,18,24\}=O_{BC}.
\end{align*}

\section{The four designs of \texorpdfstring{$S(3,4,14)$}{S(3,4,14)}}\label{app:4designs}
We list all the 91 blocks of the four designs $\BB_1,\BB_2,\BB_3$, and $\BB_4$ of
$S(3,4,14)$. To compress the description, we write $a,b,c,d,e$ for $10,11,12,13,14$,
respectively, and we write, for example, $12cd$ for $\{1,2,12,13\}$. Each list has $91$
distinct blocks, and direct enumeration shows that each of the $\binom{14}{3}=364$ triples
occurs in exactly one block. To verify that the four lists are pairwise nonisomorphic, let
$N_3(\BB)$ count the $6$-subsets of $[14]$ that contain exactly three blocks of $\BB$.
Enumeration of the $\binom{14}{6}=3003$ subsets gives
\[
N_3(\BB_1)=62,\quad N_3(\BB_2)=119,\quad
N_3(\BB_3)=51,\quad N_3(\BB_4)=98.
\]
These distinct isomorphism invariants, together with the classification in \cite{MH}, show
that the lists represent all four isomorphism classes.

\medskip
\begingroup\footnotesize
\noindent$\mathcal B_1$.
\begin{verbatim}
1236 1247 125e 128b 129a 12cd 1345 137d 138a 139c 13be 146e 148c
149b 14ad 156a 157c 1589 15bd 1679 168d 16bc 178e 17ab 19de 1ace
2348 2359 237c 23ae 23bd 245c 246d 249e 24ab 256b 2578 25ad 267a
268e 269c 279b 27de 289d 28ac 2bce 3469 347a 34bc 34de 356d 357b
358e 35ac 3678 36ab 36ce 379e 389b 38cd 39ad 4567 458d 459a 45be
468b 46ac 4789 47bd 47ce 48ae 49cd 568c 569e 579d 57ae 58ab 59bc
5cde 67be 67cd 689a 69bd 6ade 78ad 78bc 79ac 89ce 8bde 9abe abcd
\end{verbatim}
\noindent$\mathcal B_2$.
\begin{verbatim}
1236 1247 125e 128b 129a 12cd 1345 137d 138a 139c 13be 146e 148c
149b 14ad 156a 157c 1589 15bd 1679 168d 16bc 178e 17ab 19de 1ace
2348 235c 2379 23ab 23de 245b 246c 249d 24ae 2569 257a 258d 2678
26ad 26be 27bd 27ce 289e 28ac 29bc 3467 349e 34ac 34bd 356e 3578
359b 35ad 368b 369a 36cd 37ae 37bc 389d 38ce 456d 457e 458a 459c
4689 46ab 478b 479a 47cd 48de 4bce 567b 568c 579d 58be 59ae 5abc
5cde 67ac 67de 68ae 69bd 69ce 789c 78ad 79be 89ab 8bcd 9acd abde
\end{verbatim}
\noindent$\mathcal B_3$.
\begin{verbatim}
1236 1247 125e 128b 129a 12cd 1345 137d 138a 139c 13be 146e 148c
149b 14ad 156a 157c 1589 15bd 1679 168d 16bc 178e 17ab 19de 1ace
234a 235d 2379 238e 23bc 245b 2468 249d 24ce 2569 2578 25ac 267c
26ab 26de 27ae 27bd 289c 28ad 29be 346c 347b 3489 34de 356b 357a
358c 359e 3678 369d 36ae 37ce 38bd 39ab 3acd 4567 458d 459c 45ae
469a 46bd 478a 479e 47cd 48be 4abc 568e 56cd 579b 57de 58ab 59ad
5bce 67ad 67be 689b 68ac 69ce 789d 78bc 79ac 89ae 8cde 9bcd abde
\end{verbatim}
\noindent$\mathcal B_4$.
\begin{verbatim}
1236 1247 125e 128b 129a 12cd 1345 137d 138a 139c 13be 146e 148c
149b 14ad 156a 157c 1589 15bd 1679 168d 16bc 178e 17ab 19de 1ace
234d 235a 237c 238e 239b 2458 2469 24ac 24be 2567 259d 25bc 268c
26ab 26de 2789 27ae 27bd 28ad 29ce 346b 347a 3489 34ce 3569 357b
358c 35de 3678 36ae 36cd 379e 38bd 39ad 3abc 456d 457e 459c 45ab
467c 468a 478b 479d 48de 49ae 4bcd 568b 56ce 578d 579a 58ae 59be
5acd 67ad 67be 689e 69ac 69bd 78ac 79bc 7cde 89ab 89cd 8bce abde
\end{verbatim}
\endgroup


\begin{thebibliography}{99}
\bibitem{B}
D.~J. Bergstrand, New uniqueness proofs for the $(5,\,8,\,24),$ $(5,\,6,\,12)$\ and
related Steiner systems, J. Combin. Theory Ser. A {\bf 33} (1982), no.~3, 247--272;
MR0676745

\bibitem{BGM}
J. Brakensiek, S. Gopi and V. Makam, Generic Reed-Solomon codes achieve list-decoding
capacity, in {\it STOC'23---Proceedings of the 55th Annual ACM Symposium on Theory of
Computing}, 1488--1501, ACM, New York, 2023; MR4617483

\bibitem{BGM2}
J.~Brakensiek, S.~Gopi and V.~Makam,
\emph{Generic Reed-Solomon Codes Achieve List-decoding Capacity},
\href{https://arxiv.org/abs/2206.05256v4}{arXiv:2206.05256v4}, 2022.

\bibitem{CT}
T.~M. Cover and J.~A. Thomas, {\it Elements of information theory}, second edition,
Wiley-Intersci., Hoboken, NJ, 2006; MR2239987

\bibitem{Cs}
L.~Csirmaz, Exploring the entropic region,
\href{https://arxiv.org/abs/2509.12439v2}{arXiv:2509.12439v2}, 2025.

\bibitem{DEF}
M.-M. Deza, P. Erd\H os and P. Frankl, Intersection properties of systems of finite sets,
Proc. London Math. Soc. (3) {\bf 36} (1978), no.~2, 369--384; MR0476536

\bibitem{DFZ}
R.~Dougherty, C.~Freiling, and K.~Zeger,
\emph{Non-Shannon Information Inequalities in Four Random Variables},
\href{https://arxiv.org/abs/1104.3602v1}{arXiv:1104.3602v1}, 2011.

\bibitem{FOT}
P. Frankl, K. Ota and N. Tokushige, Exponents of uniform $L$-systems, J. Combin. Theory
Ser. A {\bf 75} (1996), no.~1, 23--43; MR1395757

\bibitem{FT}
P. Frankl and N. Tokushige, {\it Extremal problems for finite sets}, Student Mathematical
Library, 86, Amer. Math. Soc., Providence, RI, 2018; MR3822342

\bibitem{Furedi}
Z. F\"uredi, On finite set-systems whose every intersection is a kernel of a star,
Discrete Math. {\bf 47} (1983), no.~1, 129--132; MR0720617

\bibitem{KT}
A.~Kupavskii, N.~ Terekhov. Almost entropic matroids and $(k,L)$-systems.
\href{https://arxiv.org/abs/2610.27892v1}{arXiv:2610.27892v1}, 2026.

\bibitem{MW-S}
F.~J. MacWilliams and N.~J.~A. Sloane, {\it The theory of error-correcting codes.},
North-Holland Mathematical Library, Vol. 16, North-Holland, Amsterdam-New York-Oxford,
1977; MR0465509, MR0465510

\bibitem{MH}
N.~S. Mendelsohn and S.~H.~Y. Hung, On the Steiner systems $S(3,\,4,\,14)$ and
$S(4,\,5,\,15)$, Utilitas Math. {\bf 1} (1972), 5--95; MR0302474

\bibitem{RT}
V. R\"odl and E. Tengan, A note on a conjecture by F\"uredi, J. Combin. Theory Ser. A {\bf
113} (2006), no.~6, 1214--1218; MR2244143

\bibitem{T}
N. Tokushige, An $L$-system on the small Witt design, J. Combin. Theory Ser. A {\bf 113}
(2006), no.~3, 420--434; MR2209702

\bibitem{Y}
R.~W. Yeung, {\it A first course in information theory}, Information Technology:
Transmission, Processing and Storage, Kluwer/Plenum, New York, 2002; MR2042182

\bibitem{ZY}
Z. Zhang and R.~W. Yeung, On characterization of entropy function via information
inequalities, IEEE Trans. Inform. Theory {\bf 44} (1998), no.~4, 1440--1452; MR1665794
\end{thebibliography}
\end{document}